# A reliable numerical method for solving a certain class of singular initial value problems using reproducing kernel algorithm

Qasem Al-Haj Abdullah[1], Mohammed Al-Smadi[2], Radwan Abu-Gdairi[1], Abdel Karim Baareh[2], Asad Freihat[2] and Omar Abu Arqub[3]

[1]Department of Mathematics, Faculty of Science, Zarqa University, Zarqa 13110, Jordan
[2]Department of Applied Science, Ajloun College, Al-Balqa Applied University, Ajloun 26816, Jordan
[3]Department of Mathematics, Al-Balqa Applied University, Salt 19117, Jordan

**Abstract.** The aim of this study is to present a good modernistic strategy for solving some well-known classes of Lane-Emden type singular differential equations. The proposed approach is based on the reproducing kernel Hilbert space (RKHS) and introducing the reproducing kernel properties in which the initial conditions of the problem are satisfied. The analytical solution that obtained involves in the form of a convergent series with easily computable terms in its reproducing kernel space. The approximation solution is expressed by n-term summation of reproducing kernel functions and it is converge to the analytical solution. Our investigations indicate that there is excellent agreement between the numerical results and the RKHS method, which is applied to some examples to illustrate the accuracy, efficiency, and applicability of the method. The present work shows the potential of the RKHS technique in solving such nonlinear singular initial value problems.

## 1. Introduction

Mathematical modelling is a major challenge for contemporary scientists. Singular initial value problems (SVIPs) are one of the most powerful models arising in astrophysics, mechanics and mathematical sciences for handling Lane-Emden type equations, which are well-known singular differential equations and was studied for the first time by the astrophysicists Jonathan Lane and Robert Emden. Lane and Emden considered the thermal behavior of a spherical cloud of gas acting under the mutual attraction of its molecules and subject to the classical laws of thermodynamics. Anyhow, these types of problems of great importance in many branches of physics and engineering field. Therefore, they received special attention of scientists and researchers [1-7].

The purpose of this paper is to present an effective numerical method based on the RKHS algorithm for solving Lane-Emden singular differential equations in the following form

$$u''(x) + \frac{k}{x}u'(x) = F(x, u(x)), a < x < T, \quad (1)$$

subject to the initial conditions

$$u(a) = \alpha, u'(a) = \beta, \quad (2)$$

where $k, a, \alpha$ and $\beta$ are real finite constants, $u(x) \in W_2^3[a,T]$ is an unknown function to be determined, $F(x,y)$ is continuous term in $W_2^1[a,T]$ as $y = y(x) \in W_2^3[a,T]$, $a \le x \le T$, $\infty < y < \infty$, which is depending on the problem discussed, $W_2^1[a,T]$ and $W_2^3[a,T]$ are two reproducing kernel spaces. Throughout this analysis, we assume that $F$ and $u$ are analytic functions on the given interval as well as satisfy all the necessary requirements for the existence of a unique solution.

During the last few decades, many analytic and numeric methods were developed to study and to obtain approximate solutions for different types of Lane-Emden equations and their coupled system. The Adomian decomposition method (ADM) [8], the homotopy analysis method (HAM) [9], the variational iteration method (VIM) [10], the homotopy perturbation method (HPM) [11], the collocation method [12], and Legendre method [13] are some examples of these methods. On the other hand, other categories based upon numerical approaches for different problems can be found in [14-26] and references therein.

Reproducing kernel theory has important application in numerical analysis, computational mathematics, finance, probability and statistics [27-30]. The RKHS method was successfully used by many authors to investigate several scientific applications side by side with their theories. The present method has the following characteristics; Firstly, it is of global nature in terms of the solutions obtained as well as its ability to solve other mathematical, physical, and engineering problems; Secondly, it is accurate and the global approximation and its derivatives can be established on the whole solution domain; Finally, the method does not require discretization of the variables, and it is not affected by the computation round off errors and one is not faced with necessity of large computer memory and time. The reader is requested to refer to [31-39] in order to know more details and descriptions about the methodology of the RKHS algorithm including their history, theory, modification, characteristics and scientific applications.

This paper is organized as follows. In Section 2, we present construction of the method in the reproducing kernel space. The analytical and approximate solutions under the assumption that the solution of SIVP (1) and (2) is unique are introduced in Section 3. Implementations of the method for obtaining the solution of both linear and nonlinear Lane-Emden singular equation is presented in Section 4 upon the Hilbert space $W_2^3[a,T]$. The numerical results are reported to demonstrate the superiority and capability of the proposed scheme by considering some examples in Section 5. The last section is a brief conclusion.

## 2. Reproducing kernel space

The material in this section is basic in some sense. For the reader's convenience, it is necessary to present an appropriate brief introduction to preliminary topics from the reproducing kernel theory.

**Definition .1** [40] Let $\mathcal{H}$ be a Hilbert space of function $\emptyset: \Omega \to \mathcal{H}$ on a set $\Omega$. A function $K: \Omega \times \Omega \to \mathbb{C}$ is a reproducing kernel of the Hilbert space $\mathcal{H}$ if the following conditions are satisfied. Firstly, $K(\cdot, x) \in \mathcal{H}$ for each $x \in \Omega$. Secondly, $\langle \varphi, K(\cdot, x) \rangle = \varphi(x)$ for each $x \in \Omega$ and $\varphi \in \mathcal{H}$.

**Definition .2** [41] The inner product space $W_2^3[a,T]$ is defined as $W_2^3[a,b] = \{u: u''(x)$ are absolutely continuous functions on $[a,T], u''' \in L^2[a,b], u(a) = u'(a) = 0\}$. On the other hand, the inner product and norm in $W_2^3[a,T]$ are defined, respectively, by

$$\langle u(x), v(x) \rangle_{W_2^3} = \sum_{i=0}^{2} u^{(i)}(a) v^{(i)}(a) + \int_0^1 u'''(s) v'''(s) ds, \tag{3}$$

and $\|u\|_{W_2^3} = \sqrt{\langle u, u \rangle_{W_2^3}}$, where $u, v \in W_2^3[a,T]$.

The Hilbert space $W_2^3[a,T]$ is called a reproducing kernel if for each fixed $x \in [a,T]$, there exist $R(x,y) \in W_2^3[a,T]$ (simply $R_x(y)$) such that $\langle u(y), R_x(y)\rangle_{W_2^3} = u(x)$ for any $u(y) \in W_2^3[a,T]$ and $y \in [a,T]$. Next theorem utilize the reproducing kernel function $R_x(y)$ on the space $W_2^3[a,T]$.

**Theorem .1** The Hilbert space $W_2^3[a,T]$ is a complete reproducing kernel and its reproducing kernel function $R_x(y)$ can be written as

$$R_x(y) = \begin{cases} a_1(x) + a_2(x)y + a_3(x)y^2 + a_4(x)y^3 + a_5(x)y^4 + a_6(x)y^5, & y \leq x, \\ b_1(x) + b_2(x)y + b_3(x)y^2 + b_4(x)y^3 + b_5(x)y^4 + b(x)y^5, & y > x, \end{cases} \quad (4)$$

where $a_i(x)$ and $b_i(x)$, $i = 1,2,\ldots,6$, are unknown coefficients of $R_x(y)$ and will be given by the following assumptions:

**Proof.** The proof of the completeness and reproducing property of $W_2^3[a,T]$ is similar to the proof in [39]. Now, let's assume that $R_x(y) \in W_2^3[a,T]$ satisfies the generalized differential equations

$$\begin{cases} \partial_y^6 R_x(y) = -\delta(x-y), R_x(a) = R_x(T) = 0, \\ \partial_y^2 R_x(a) - \partial_y^3 R_x(a) = 0, \\ \partial_y^5 R_x(T) = \partial_y^4 R_x(T) = \partial_y^3 R_x(T) = 0. \end{cases} \quad (5)$$

where $\delta$ is the Dirac delta function.

On the other hand, for $x \neq y$, $R_x(y)$ is the solution of the constant differential equation $\partial_y^6 R_x(y) = 0$ subject to the conditions (5). That is, the characteristic equation is given by $\lambda^6 = 0$ and the eigenvalues are $\lambda = 0$ with multiplicity 6. Hence, the general solution can be written as in Eq. (4). But on the other aspect as well, let $R_x(y)$ satisfy the equation $\partial_y^m R_x(x+0) = \partial_y^m R_x(x-0)$, $m = 0,1,2,3,4$. Integrating $\partial_y^6 R_x(y) = -\delta(x-y)$ from $x-\varepsilon$ to $x+\varepsilon$ with respect to $y$ and let $\varepsilon \to 0$, we have the jump degree of $\partial_y^5 R_x(y)$ at $y = x$ given by $\partial_y^5 R_x(x+0) - \partial_y^5 R_x(x-0) = -1$. Hence, through these descriptions the unknown coefficients $a_i(x)$ and $b_i(x)$, $i = 1,2,\ldots,6$, of Eq. (4) can be obtained. This completes the proof.

**Remark .1** The unique representation of the reproducing kernel $R_x(y)$ is provided

$$R_x(y) = \begin{cases} \dfrac{-1}{120}(a-y)^2 \begin{pmatrix} 6a^3 + 5xy^2 - y^3 - 10x^2(3+y) - 3a^2( \\ 10 + 5x + y) + 2a(5x^2 - y^2 + 5x(6+y)) \end{pmatrix}, & y \leq x, \\ \dfrac{1}{120}(a-x)^2 \begin{pmatrix} -6a^3 + x^3 - 5xy^2 + 30y^2 + 10x^2y + 3a^2(x) \\ +10 + 5y + 2a(x^2 - 5x^2 - y^2 - 5y(6+y)) \end{pmatrix}, & y > x. \end{cases}$$

The reproducing kernel function $R_x(y)$ possess some important properties such as: is symmetric, unique and nonnegative for any fixed $x \in [a,T]$.

**Definition 3:** [42] The inner product space $W_2^1[a,T]$ is defined as $W_2^1[a,T] = \{u(x): u'(x) \text{ is absolutely continuous real-valued function, } u'(x) \in L^2[a,T]\}$. On the other hand, the inner product and norm in $W_2^1[a,T]$ are defined, respectively, by

$$\langle u(x), v(x)\rangle_{W_2^1} = u(a)v(a) + \int_a^T u'(t)v'(t)dt,$$

and $\|u\| = \sqrt{\langle u(x), u(x) \rangle_{W_2^1}}$, where $u, v \in W_2^1[a, T]$.

**Remark .2** In [42], it has been proved that the Hilbert space $W_2^1[a, T]$ is a complete reproducing kernel and its reproducing kernel function is given by

$$G_x(y) = \frac{1}{2\sinh(T-a)}\left[\cosh(x+y-(a+T)) + \cosh(|x-y|+a-T)\right].$$

## 3. Formulation of analytical solution

In this section, the formulation of exact and approximate solutions of Eqs. (1) and (2) together with implementation method are given in the reproducing kernel space $W_2^3[a, T]$. After that, we construct an orthogonal function system of the space $W_2^3[a, T]$ based upon the use of the Gram-Schmidt orthogonalization process.

To do this, we define a differential operator $L$ as

$$L: W_2^3[a, T] \to W_2^1[a, T]$$

such that

$$Lu(x) = u''(x) + \frac{k}{x}u'(x).$$

After homogenization of the initial conditions, the Eqs. (1) and (2) can be converted into the equivalent form as follows:

$$Lu(x) = F(x, u(x)); \quad a < x < T,$$
$$u(a) = 0, u'(a) = 0,$$
(6)

where $u(x) \in W_2^3[a, T]$ and $F(x, y) \in W_2^1[a, T]$ for $y = y(x) \in W_2^3[a, T]$, $\infty < y < \infty$. It is easy to show that $L$ is a bounded linear operator from the space $W_2^3[a, T]$ into space $W_2^1[a, T]$.

Initially, we construct an orthogonal function system of $W_2^3[a, T]$. To do so, put $\varphi_i(x) = G_{x_i}(x)$ and $\psi_i(x) = L^*\varphi_i(x)$, where $\{x_i\}_{i=1}^\infty$ is dense on $[a, T]$ and $L^*$ is the adjoint operator of $L$. In terms of the properties of reproducing kernel $G_x(y)$, one obtains $\langle u(x), \psi_i(x) \rangle_{W_2^3} = \langle u(x), L^*\varphi_i(x) \rangle_{W_2^3} = \langle Lu(x), \varphi_i(x) \rangle_{W_2^1} = Lu(x_i)$, $i = 1, 2, \ldots$ For the orthonormal function system $\{\bar{\psi}_i(x)\}_{i=1}^\infty$ of the space $W_2^3[a, T]$, it can be derived from the Gram-Schmidt orthogonalization process of $\{\psi_i(x)\}_{i=1}^\infty$ as follows:

$$\bar{\psi}_i(x) = \sum_{k=1}^{i} \beta_{ik} \psi_k(x), \tag{7}$$

where $\beta_{ik}$ are orthogonalization coefficients ($\beta_{ii} > 0, i = 1, 2, \ldots, n$), which are given by the following subroutine:

$\beta_{11} = \frac{1}{\|\psi_i\|}$; $\beta_{ij} = \frac{1}{\|\psi_i\|}$ for $i = j = 1$; $\beta_{ij} = \frac{1}{d_{ik}}$ for $i = j \neq 1$; $\beta_{ij} = -\frac{1}{d_{ik}}\sum_{k=j}^{i-1} c_{ik}\beta_{jk}$ $i > j$,

such that $d_{ik} = \sqrt{\|\psi_i\|^2 - \sum_{k=1}^{i-1} c_{ik}^2}$, $c_{ik} = \langle \psi_i, \bar{\psi}_k \rangle_{W_2^3}$ and $\{\psi_i(x)\}_{i=1}^\infty$ is the orthonormal system in the space $W_2^3[a, T]$.

It is easy to see that, $\psi_i(x) = L^*\varphi_i(x) = \langle L^*\varphi_i(x), R_x(y) \rangle_{W_2^3} = \langle \phi_i(x), L_y R_x(y)\rangle_{W_2^1} = L_y R_x(y)\big|_{y=x_i} \in W_2^3[a,T]$. Thus, $\psi_i(x)$ can be written in the form $\psi_i(x) = L_y R_x(y)\big|_{y=x_i}$, where $L_y$ indicates that the operator $L$ applies to the function of $y$.

**Theorem .2** If $\{x_i\}_{i=1}^{\infty}$ is dense on $[a,T]$, then $\{\psi_i(x)\}_{i=1}^{\infty}$ is a complete function system of the space $W_2^3[a,T]$.

**Proof.** For each fixed $u(x) \in W_2^3[a,T]$, let $\langle u(x), \psi_i(x)\rangle_{W_2^3} = 0$, $i = 1,2,\ldots$. In other word, one has, $\langle u(x), \psi_i(x)\rangle_{W_2^3} = \langle Lu(x), \phi_i(x)\rangle_{W_2^1} = Lu(x_i) = 0$. Note that $\{x_i\}_{i=1}^{\infty}$ is dense on $[a,T]$, therefore $Lu(x) = 0$. It follows that $u(x) = 0$ from the existence of $L^{-1}$. So, the proof of the theorem is complete.

**Lemma .1** If $u(x) \in W_2^3[a,T]$, then there exists positive constants $M$ such that $\|u^{(i)}(x)\|_c \leq M\|u(x)\|_{W_2^3}$, $i = 0,1$, where $\|u(x)\|_c = \max_{a\leq x\leq T}|u(x)|$.

**Proof.** For any $x, y \in [a,T]$, we have $u^{(i)}(x_1) = \langle u(y), \partial_x^i K(x,y)\rangle_{W_2^3}, i = 0,1$. By the expression form of $R_x(y)$, it follows that $\|\partial_x^i R_x(y)\|_{W_2^3} \leq M_i, i = 0,1$. Thus, $|u^{(i)}(x)| = \left|\langle u(x), \partial_x^i R_x(x)\rangle_{W_2^3}\right| \leq \|\partial_x^i R_x(x)\|_{W_2^3}\|u(x)\|_{W_2^3} \leq M_i\|u(x)\|_{W_2^2}, i = 0,1$.

The internal structure of the following theorem is as follows: firstly, we will give the representation of the exact solution of Eqs. (1) and (2) in the space $W_2^3[a,T]$. After that, the convergence of approximate solution $u_n(x)$ to the analytic solution $u(x)$ will be proved.

**Theorem .3** For each $u(x)$ in $W_2^3[a,T]$, the series $\sum_{i=1}^{\infty}\langle u(x), \bar{\psi}_i(x)\rangle \bar{\psi}_i(x)$ is convergent in the sense of the norm of $W_2^3[a,T]$. On the other hand, if $\{x_i\}_{i=1}^{\infty}$ is dense on $[a,T]$, then the following are hold:

(i) the exact solution of Eq. (6) could be represented by

$$u(x) = \sum_{i=1}^{\infty}\sum_{k=1}^{i} \beta_{ik} F(x_k, u(x_k))\bar{\psi}_i(x), \qquad (8)$$

(ii) the approximate solution of Eq. (6) is also represented by

$$u_n(x) = \sum_{i=1}^{n}\sum_{k=1}^{i} \beta_{ik} F(x_k, u(x_k))\bar{\psi}_i(x), \qquad (9)$$

and $u_n^{(i)}(x), i = 0,1,2$ are converging uniformly to the exact solution $u(x)$ and all its derivative as $n \to \infty$, respectively.

**Proof.** For the first part, let $u(x)$ be solution of Eq. (6) in the space $W_2^3[a,T]$. Since $u(x) \in W_2^3[a,T]$, then $u(x) = \sum_{i=1}^{\infty}\langle u(x), \bar{\psi}_i(x)\rangle \bar{\psi}_i(x)$ is the Fourier series expansion about normal orthogonal system $\{\psi_i(x)\}_{i=1}^{\infty}$, but $W_2^3[a,T]$ is the Hilbert space, then the series $\sum_{i=1}^{\infty}\langle u(x), \bar{\psi}_i(x)\rangle \bar{\psi}_i(x)$ is convergent in the sense of $\|.\|_{W_2^3}$. On the other hand, using Eq. (7), it easy to see that

$$u(x) = \sum_{i=1}^{\infty}\langle u(x), \bar{\psi}_i(x)\rangle_{W_2^3}\bar{\psi}_i(x) = \sum_{i=1}^{\infty}\langle u(x), \sum_{k=1}^{i}\beta_{ik}\psi_k(x)\rangle_{W_2^3}\bar{\psi}_i(x)$$

$$= \sum_{i=1}^{\infty}\sum_{k=1}^{i} \beta_{ik} \langle u(x), \psi_k(x)\rangle_{W_2^3} \bar{\psi}_i(x) = \sum_{i=1}^{\infty}\sum_{k=1}^{i} \beta_{ik} \langle u(x), L^*\varphi_i(x)\rangle_{W_2^3} \bar{\psi}_i(x)$$

$$= \sum_{i=1}^{\infty}\sum_{k=1}^{i} \beta_{ik} \langle Lu(x), \varphi_k(x)\rangle_{W_2^3} \bar{\psi}_i(x) = \sum_{i=1}^{\infty}\sum_{k=1}^{i} \beta_{ik} \langle F(x_k, u(x_k)), \varphi_k(x)\rangle_{W_2^3} \bar{\psi}_i(x)$$

$$= \sum_{i=1}^{\infty}\sum_{k=1}^{i} \beta_{ik} F(x_k, u(x_k)) \bar{\psi}_i(x).$$

Therefore, the form of Eq. (8) is the exact solution of Eq. (6). For the second part, it easy to see that by Lemma 1, for any $x \in [a, T]$, we have

$$|u_{s,n}(x) - u_s(x))| = |\langle u_{s,n}(x) - u_s(x), R_x(x)\rangle_{W_2^3}| \leq \|R_x(x)\|_{W_2^3} \|u_{s,n}(x) - u_s(x)\|_{W_2^2}$$

$$\leq M_0^{\{s\}} \|u_{s,n}(x) - u_s(x)\|_{W_2^3} \quad s = 0,1.$$

On the other hand,

$$|u_n^{(i)}(x) - u^{(i)}(x)| = |\langle u_n(x) - u(x), R_x^{(i)}(x)\rangle_{W_2^3}| \leq \|\partial_x^i R_x(x)\|_{W_2^2} \|u_n(x) - u(x)\|_{W_2^3}$$

$$\leq M_i \|u_n(x) - u(x)\|_{W_2^3}, s = 0,1.$$

where $M_i, i = 0,1,2$, are positive constants. Hence, if $\|u_n(x) - u(x)\|_{W_2^3} \to 0$ as $n \to \infty$, the approximate solution $u_n(x)$ and $u_n^{(i)} i = 0,1,2$, are converge uniformly to the exact solution $u_n(x)$ and all its derivative, respectively. So, the proof of the theorem is complete.

**Remark 3:** We mention here that, the approximate solution $u_n(x)$ in Eq. (9) can be obtained directly by taking finitely many terms in the series representation for $u(x)$ of Eq. (8).

## 4. Implementations with convergence analysis

In this section, an iterative method of obtaining the solution of Eq. (6) is represented in the reproducing kernel space $W_2^3[a, T]$ for linear and nonlinear case. Initially, we will mention the following remark about the exact and approximate solutions of Eqs. (1) and (2).

**Remark .4** In order to apply the RKHS technique for solve Eqs. (1) and (2), we have the following two cases based on the structure of the function $F(x, u)$.

- **Case 1:** If Eq. (1) is linear, then the exact and approximate solutions can be obtained directly from Eqs. (8) and (9), respectively.
- **Case 2:** If Eq. (1) is nonlinear, then in this case the exact and approximate solutions can be obtained by using the following iterative algorithm:

**Algorithm .1** According to Eq. (8), the representation of the solution of problem (1) can be denoted by

$$u(x) = \sum_{i=1}^{\infty} B_i \bar{\Psi}_i(x), \tag{10}$$

where $B_i = \sum_{k=1}^{i} \beta_{ik} f(x_k, u_{k-1}(x_k))$. In fact, $B_i, i = 1,2, ...$, in equation (10) are unknown, so we will approximate them using the known $A_i$ as follows: For a numerical computations, let the initial function $u_0(x_1) = 0$, set $u_0(x_1) = u(x_1)$, and difne the $n$-term approximation to $y_s(x)$ by

$$u_n(x) = \sum_{i=1}^{n} A_i \overline{\Psi}_i(x), \quad (11)$$

where the coefficients $A_i$ of $\overline{\Psi}_i(x), i = 1,2, ..., n$, are given by

$$\begin{cases} A_1 = \beta_{11} f(x_1, u_0(x_1)), u_1(x) = A_1 \overline{\Psi}_1(x), \\ A_2 = \sum_{k=1}^{2} \beta_{2k} f(x_1, u_{k-1}(x_k)), u_2(x) = \sum_{i=1}^{2} A_i \overline{\Psi}_i(x), \\ \vdots \\ u_{n-1}(x) = \sum_{i=1}^{n-1} A_i \overline{\Psi}_i(x), A_n = \sum_{k=1}^{n} \beta_{nk} f(x_1, u_{k-1}(x_k)). \end{cases} \quad (12)$$

Consequently, the unknown coefficients $B_i, i = 1,2, ...$, in Eq. (10) will be approximate using the known coefficients $A_i, i = 1,2, ...$, that given in Eq. (12). However, in the iterative process of the series (11), we can guarantee that the approximation $u_n(x)$ satisfies the periodic boundary condition (2).

Now, we will proof that $u_n(x)$ in the iterative formula (11) are converge to the exact solution $u(x)$ of Eq. (1). In fact, this result is a fundamental in the RKHS theory and its applications. The next two lemmas are collected in order to prove the recent theorem.

**Lemma .2** If $\|u_n(x) - u(x)\|_{W_2^3} \to 0, x_n \to y, (n \to \infty)$, and $f(x, z)$ is continuous in $[a, T]$ with respect to $x, z$ for $x \in [a, T], z \in (-\infty, \infty)$, then the following are hold in the sense of the norm of $W_2^3[a, T]$:

(a) $u_{n-1}(x_n) \to u(y)$ as $n \to \infty$.

(b) $f(x_n, u_{n-1}(x_n)) \to f(y, u(y))$, as $n \to \infty$.

**Proof**: For part (a), note that $|u_{n-1}(x_n) - u(y)| = |u_{n-1}(x_n) - u_{n-1}(y) + u_{n-1}(y) - u(y)| \leq |u_{n-1}(x_n) - u_{n-1}(y)| + |u_{n-1}(y) - u(y)|$. By reproducing property of $R(x, y)$, we have $u_{n-1}(x_n) = \langle u_{n-1}(x), R(x_n, x) \rangle_{W_2^3}$ and $u_{n-1}(y) = \langle u_{n-1}(x), R(y, x) \rangle_{W_2^3}$. Thus, $|u_{n-1}(x_n) - u_{n-1}(y)| = \left|\langle u_{n-1}(x), R(x_n, x) - R(y, x) \rangle_{W_2^3}\right| \leq \|u_{n-1}(x)\|_{W_2^3} \|R(x_n, x) - R(y, x)\|_{W_2^3}$.

From the symmetry of $R(x, y)$, it follows that $\|R(x_n, x) - R(y, x)\|_{W_2^3} \to 0$ as $x_n \to y, n \to \infty$. Hence, $|u_{n-1}(x_n) - u_{n-1}(y)| \to 0$ as soon as $x_n \to y, (n \to \infty)$. On the other hand, for any $x \in [0,1]$, by using Theorem 2, it holds that $|u_{n-1}(y) - u(y)| \to 0$ as $n \to \infty$. Therefore, $u_{n-1}(x_n) \to u(y)$ in the sense of $\|\cdot\|_{W_2^3}$ as $x_n \to y$ and $n \to \infty$. Thus, for part (b), by means of the continuation of $f(\cdot)$, it is obtained that $f(x_n, u_{n-1}(x_n)) \to f(y, u(y))$ as $x_n \to y$ and $n \to \infty$.

**Lemma .3** $Lu_n(x_j) = Lu(x_j) = f\left(x_j, u_{j-1}(x_j)\right), j \leq n$, holds.

**Proof**: The proof of $u_n(x_j) = f\left(x_j, u_{j-1}(x_j)\right)$ will be obtained by induction as follows: if $j \leq n$, then

$$Lu_n(x_j) = \sum_{i=1}^{n} A_i L \overline{\Psi}_i(x) = \sum_{i=1}^{n} A_i \langle L \overline{\Psi}_i(x), \Phi_j(x) \rangle_{W_2^1} = \sum_{i=1}^{n} A_i \langle \overline{\Psi}_i(x), L^* \Phi_j(x) \rangle_{W_2^3}$$

$$= \sum_{i=1}^{n} A_i \langle \overline{\Psi}_i(x), \Psi_j(x) \rangle_{W_2^3}.$$

That is, $Lu_n(x_j) = \sum_{i=1}^{n} A_i \langle \overline{\Psi}_i(x), \Psi_j(x) \rangle_{W_2^3}$. Multiplying both sides of the previous relation by $\beta_{jl}$, summing for $l$ from 1 to $j$, and using the orthogonality of $\{\overline{\Psi}_i(x)\}_{i=1}^{\infty}$, yields that

$$\sum_{l=1}^{j} \beta_{jl} Lu_n(x_l) = \sum_{i=1}^{n} A_i \langle \overline{\Psi}_i(x), \sum_{l=1}^{j} \beta_{jl} \Psi_l(x) \rangle_{W_2^3} = \sum_{i=1}^{\infty} A_i \langle \overline{\Psi}_i(x), \overline{\Psi}_j(x) \rangle_{W_2^3} = A_j$$

$$= \sum_{l=1}^{j} \beta_{jl} f(x_l, u_{l-1}(x_l)).$$

If $j = 1$, then $Lu_n(x_1) = f(x_1, u_0(x_1))$. Besides, if $j = 2$, then $\beta_{21} Lu_n(x_1) + \beta_{22} Lu_n(x_2) = \beta_{21} f(x_1, u_0(x_1)) + \beta_{22} f(x_2, u_1(x_2))$, that is, $Lu_n(x_2) = f(x_2, u_1(x_2))$. By the same manner, yields that $Lu_n(x_j) = f(x_j, u_{j-1}(x_j))$ for $j \leq n$.

But on the other aspects as well, from Theorem 2 as well by taking limits in Eq. (11), we have $u(x) = \sum_{i=1}^{\infty} A_i \overline{\Psi}_i(x)$. Thus, $u_n(x) = P_n u(x)$, where $P_n$ is an orthogonal projector from the space $W_2^3[a, T]$ to $Span\{\Psi_1, \Psi_2, \ldots, \Psi_n\}$. Therefore,

$$Lu_n(x_j) = \langle Lu_n(x), \Phi_j(x) \rangle_{W_2^1} = \langle u_n(x), L^*\Phi_j(x) \rangle_{W_2^3} = \langle P_n u(x), \Psi_j(x) \rangle_{W_2^3} = \langle u(x), P_n \Psi_j(x) \rangle_{W_2^3}$$
$$= \langle u(x), \Psi_j(x) \rangle_{W_2^3} = \langle u(x), L^*\Phi_j(x) \rangle_{W_2^3} = \langle Lu(x), \Phi_j(x) \rangle_{W_2^1} = Lu(x_j).$$

**Theorem .4** If $\|u_n(x)\|_{W_2^3}$ is bounded and $\{x_i\}_{i=1}^{\infty}$ is dense on $[a, T]$, then the $n$-term approximate solution $u_n(x)$ in the iterative formula (11) is convergent to the exact solution $u(x)$ of Eq. (6) in the space $W_2^3[a, T]$ and $u(x) = \sum_{i=1}^{\infty} A_i \overline{\Psi}_i(x)$, where $A_i, i = 1, 2, \ldots$ are given by Equation (12).

***Proof***: The proof consists of the following three steps: firstly, we will prove the convergence of $\{u_n\}_{n=1}^{\infty}$ in Eq. (11) is monotone increasing in the sense of $\|\cdot\|_{W_2^3}$. By Theorem 2, $\{\overline{\Psi}_i(x)\}_{i=1}^{\infty}$ is the complete orthonormal system in the space $W_2^3[a, T]$. Hence, we have $\|u_{n+1}\|_{W_2^3}^2 = \|u_n\|_{W_2^3}^2 + (A_{n+1})^2 = \|u_{n-1}\|_{W_2^3}^2 + (A_n)^2 + (A_{n+1})^2 = \cdots = \|u_0\|_{W_2^3}^2 + \sum_{i=1}^{n+1}(A_i)^2$. Therefore, the sequence $\|u_n\|_{W_2^3}$ is monotone increasing and from the boundedness of $\|u_n\|_{W_2^3}$, we have $\sum_{i=1}^{\infty}(A_i)^2 < \infty$, that is, $\{A_i\}_{i=1}^{\infty} \in l^2 (i = 1, 2, \ldots)$. Hence, $\|u_n\|_{W_2^3}$ is convergent as $n \to \infty$. Let $m > n$, for $(u_m - u_{m-1}) \perp (u_{m-1} - u_{m-2}) \perp \cdots \perp (u_{n+1} - u_n)$, it follows that

$$\|u_m(x) - u_n(x)\|_{W_2^3}^2 = \|u_m(x) - u_{m-1}(x) + u_{m-1}(x) - \cdots + u_{n+1}(x) - u_n(x)\|_{W_2^3}^2$$

$$\leq \|u_m(x) - u_{m-1}(x)\|_{W_2^3}^2 + \cdots + \|u_{n+1}(x) - u_n(x)\|_{W_2^3}^2 = \sum_{i=n+1}^{m}(A_i)^2 \to 0, (n \to \infty).$$

Considering the completeness of $W_2^3[a, T]$, there exists $u(x) \in W_2^3[a, T]$ such that $u_n(x) \to u(x)$ as $n \to \infty$ in sense of $\|\cdot\|_{W_2^3}$. Now, we will prove that $u(x)$ is the solution of Eq. (6). Since $\{x_i\}_{i=1}^{\infty}$ is dense on compact interval $[a, T]$, thus for any $x \in [a, T]$, there exists subsequence $\{x_{n_j}\}$ such that $x_{n_j} \to x$, as $j \to \infty$. From Lemma 3, $Lu_n(x_{n_j}) = f(x_{n_j}, u_{j-1}(x_{n_j}))$. Hence, let $j \to \infty$, we have $Lu(x) = f(x, u(x))$. That is, $u(x)$ is solution of Equation (9). The proof is complete.

It obvious that, if we let $u(x)$ denote the exact solution of Eq. (6), $u_n(x)$ denote the approximate solution obtained by the RKHS method as given by Eq. (11), and $r_n(x)$ is the difference between $u_n(x)$ and $u(x)$, where $x \in [a, T]$, then $\|r_n(x)\|_{W_2^3}^2 = \|u(x) - u_n(x)\|_{W_2^3}^2 = \|\sum_{i=n+1}^{\infty} A_i \overline{\Psi}_i(x)\|_{W_2^3}^2 =$

$\sum_{i=n+1}^{\infty}(A_i)^2$ and $\|r_{n-1}(x)\|_{W_2^3}^2 = \sum_{i=n}^{\infty}(A_i)^2$. Thus, $\|r_n(x)\|_{W_2^3} \le \|r_{n-1}(x)\|_{W_2^3}$. Consequently, the error $r_n(x)$ is monotone decreasing in the sense of $\|.\|_{W_2^3}$.

## 5. Applications and numerical results

In this final section, we consider the following three examples in order to illustrate the performance of RKHS method in approximating the solution for linear and nonlinear Lane-Emden SIVPs and justify the accuracy and applicability of the method. On the other hand, results obtained by the method are compared with the exact solution of each example by computing the exact error and are found to be in good agreement with each other. In the process of computation, all the symbolic and numerical computations performed by using Mathematica 7.0 software package.

**Example 1.** Consider the linear Lane-Emden equation

$$u''(x) + \frac{2}{x}u'(x) + u(x) = x^4 + x^3 + 12x + 6, 0 \le x \le 1,$$

subject to the initial conditions

$$u(0) = 0, u'(0) = 0.$$

The exact solution is $u(x) = x^3 + x^2$.

Using RKHS method and taking $x_i = (i/n), i = 1,2,\ldots,n$ $(n = 100)$ with the reproducing kernel function $R_x(y)$ on $[0,1]$ so that

$$R_x(y) = \begin{cases} \frac{1}{120}y^2(y^3 - 5xy^2 + 10x^2(3 + y)), y \le x, \\ \frac{1}{120}x^2(x^3 - 5xy^2 + 30y^2 + 10x^2y), y > x. \end{cases}$$

The approximate solution $u_n(x)$ is calculated by Eq. (9), as well the numerical results at some selected grid points are given in Table 1.

**Table 1:** Numerical results for Example 1.

| $x_i$ | Exact solution | Approximate solution | Absolute Error | Relative error |
|---|---|---|---|---|
| 0.16 | 0.029696 | 0.0296965312802189 | $5.31280 \times 10^{-7}$ | $1.78906 \times 10^{-5}$ |
| 0.32 | 0.135168 | 0.1351686285866068 | $6.28587 \times 10^{-7}$ | $4.65041 \times 10^{-6}$ |
| 0.48 | 0.340992 | 0.3409929556422331 | $9.55642 \times 10^{-7}$ | $2.80254 \times 10^{-6}$ |
| 0.64 | 0.671744 | 0.6717457218370337 | $1.72184 \times 10^{-6}$ | $2.56323 \times 10^{-6}$ |
| 0.80 | 1.152000 | 1.1520031841576759 | $3.18416 \times 10^{-6}$ | $2.76403 \times 10^{-6}$ |
| 0.96 | 1.806336 | 1.8063416256904894 | $5.62569 \times 10^{-6}$ | $3.11442 \times 10^{-6}$ |

**Example 2.** Consider the nonlinear Lane-Emden equation

$$u''(x) + \frac{2}{x}u'(x) + 4(2e^{u(x)} + e^{\frac{1}{2}u(x)}) = 0, 0 \le x \le 1,$$

subject to the initial conditions

$$u(0) = 0, u'(0) = 0.$$

The exact solution is $u(x) = -2\ln(x^2 + 1)$.

The approximate solution $u_n(x)$ is calculated by Eq. (11), as well the numerical results at some selected grid points are given in Table 2.

Table 2: Numerical results for Example 2.

| $x_i$ | Exact solution | Approximate solution | Absolute Error | Relative error |
|---|---|---|---|---|
| 0.16 | −0.0505556 | −0.050554955996037 | $6.58372 \times 10^{-7}$ | $1.30227 \times 10^{-5}$ |
| 0.32 | −0.1949792 | −0.194977745554394 | $1.49700 \times 10^{-7}$ | $7.67774 \times 10^{-6}$ |
| 0.48 | −0.4146786 | −0.414675836611507 | $2.80293 \times 10^{-7}$ | $6.75928 \times 10^{-6}$ |
| 0.64 | −0.6866119 | −0.686608107905340 | $3.84449 \times 10^{-6}$ | $5.59922 \times 10^{-6}$ |
| 0.80 | −0.9893925 | −0.989388262187447 | $4.22148 \times 10^{-6}$ | $4.26674 \times 10^{-6}$ |
| 0.96 | −1.3063163 | −1.306312360243357 | $3.98444 \times 10^{-6}$ | $3.05014 \times 10^{-6}$ |

**Example 3.** Consider the nonlinear Lane-Emden equation

$$u''(x) + \frac{8}{x}u'(x) + +9\pi u(x) + 2\pi u(x)\ln(u(x)) = 0, 0 \leq x \leq 1,$$

subject to the initial conditions

$$u(0) = 0, u'(0) = 0.$$

The exact solution is $u(x) = e^{-\frac{\pi}{2}x^2}$.

The approximate solution $u_n(x)$ is calculated by Eq. (11), as well the numerical results at some selected grid points are given in Table 3.

Table 3: Numerical results for Example 3.

| $x_i$ | Exact solution | Approximate solution | Absolute Error | Relative error |
|---|---|---|---|---|
| 0.16 | 0.960585 | 0.960585516066161 | $1.13439 \times 10^{-7}$ | $1.18094 \times 10^{-7}$ |
| 0.32 | 0.851420 | 0.851420344457335 | $1.81599 \times 10^{-7}$ | $2.13289 \times 10^{-7}$ |
| 0.48 | 0.696344 | 0.696344424224703 | $4.14348 \times 10^{-7}$ | $5.95034 \times 10^{-7}$ |
| 0.64 | 0.525504 | 0.525504359322166 | $7.26418 \times 10^{-7}$ | $1.38233 \times 10^{-7}$ |
| 0.80 | 0.365931 | 0.365932098266620 | $7.91325 \times 10^{-7}$ | $2.16250 \times 10^{-7}$ |
| 0.96 | 0.235123 | 0.235123288472108 | $1.46755 \times 10^{-7}$ | $6.24161 \times 10^{-7}$ |

## 5. Conclusions

In this paper, we construct a reproducing kernel space in which each function satisfies the initial value conditions of considered problems. In this space, a numerical algorithm is presented for solving a class of Lane-Emden type singular differential equations. The analytical solution is given with series form in $W_2^3[a,T]$. The approximate solution obtained by present algorithm converges to analytical solution uniformly. The numerical results are displayed to demonstrate the validity of the present algorithm.